\documentclass[12pt]{article}
\usepackage{amsmath}
\usepackage{theorem,amsfonts,amssymb}

\newtheorem{theorem}{Theorem}[section]
\newtheorem{proposition}[theorem]{Proposition}
\newtheorem{lemma}[theorem]{Lemma}
\newtheorem{corollary}[theorem]{Corollary}
\theorembodyfont{\upshape}
\newtheorem{definition}{Definition}[theorem]
\newtheorem{example}[theorem]{Example}
\newtheorem{remark}[theorem]{Remark}
\newtheorem{proof}{Proof}

\newtheorem{acknowledgement}{Acknowledgement}

\newcommand{\bt}{\begin{theorem}}
\newcommand{\et}{\end{theorem}}
\newcommand{\bl}{\begin{lemma}}
\newcommand{\el}{\end{lemma}}
\newcommand{\bp}{\begin{proposition}}
\newcommand{\ep}{\end{proposition}}
\newcommand{\bd}{\begin{definition}}
\newcommand{\bex}{\begin{example}}
\newcommand{\eex}{\end{example}}
\newcommand{\ed}{\end{definition}}
\newcommand{\br}{\begin{remark}}
\newcommand{\er}{\end{remark}}
\newcommand{\bc}{\begin{corollary}}
\newcommand{\ec}{\end{corollary}}
\newcommand{\bo}{\begin{proof}}
\newcommand{\eo}{\end{proof}}
\newcommand{\be}{\begin{enumerate}}
\newcommand{\ee}{\end{enumerate}}

\newcommand{\Inn}{{\rm Inn}}

\newcommand{\Ker}{{\rm Ker}}
\newcommand{\Aut}{{\rm Aut}}

\newcommand\du{\,{\rm d}}
\newcommand{\cG}{{\cal G}}
\newcommand{\K}{{\cal K}}
\newcommand{\supp}{{\rm supp}}
\newcommand{\Z}{{\mathbb Z}}
\newcommand{\Q}{{\mathbb Q}}
\newcommand{\N}{{\mathbb N}}
\newcommand{\mF}{{\mathbb F}}

\newcommand{\C}{{\mathbb C}}
\newcommand{\T}{{\mathbb T}}

\title{Some Properties of Distal Actions on Locally Compact Groups}
\author{C.\ R.\ E.\ Raja and Riddhi Shah}
\date{ }

\begin{document}
\maketitle

\let\epsi=\epsilon
\let\vepsi=\varepsilon
\let\lam=\lambda
\let\Lam=\Lambda
\let\ap=\alpha
\let\vp=\varphi
\let\ra=\rightarrow
\let\Ra=\Rightarrow
\let\da=\downarrow
\let\Llra=\Longleftrightarrow
\let\Lla=\Longleftarrow
\let\lra=\longrightarrow
\let\Lra=\Longrightarrow
\let\ba=\beta
\let\ga=\gamma
\let\Ga=\Gamma
\let\un=\upsilon
\let\mi=\setminus
\let\ol=\overline
\let\ot=\odot

\begin{abstract} We consider the actions of (semi)groups on a locally compact group by automorphisms.
We show the equivalence of distality and pointwise distality for the actions of a certain class of groups.
We also show that a compactly generated locally compact
group of polynomial growth has a compact normal subgroup $K$ such
that $G/K$ is distal and the conjugacy action of $G$ on $K$ is
ergodic; moreover, if $G$ itself is (pointwise) distal then $G$ is
Lie projective.  We prove a decomposition theorem for contraction
groups of an automorphism under certain conditions. We give a
necessary and sufficient condition for distality of an automorphism
in terms of its contraction group.  We compare classes of
(pointwise) distal groups and groups whose closed subgroups are
unimodular.  In particular, we study relations between distality,
unimodularity and contraction subgroups.\end{abstract}



\begin{section}{Introduction}

Let $\Ga$ be a (topological) semigroup acting on a Hausdorff space
$X$ by continuous self-maps.  We say that the action of $\Ga$ on $X$
is {\it distal} if for any two distinct points $x, y \in X$, the
closure  of $\{ (\ga (x), \ga (y)) \mid \ga \in \Ga \}$ does not
intersect the diagonal $\{ (a, a ) \mid a \in X \}$ and we say that
the action of $\Ga$ on $X$ is {\it pointwise distal} if for each
$\gamma\in\Gamma$, the action of $\{ \gamma ^n\}_{n\in\N}$ on $X$ is
distal. The notion of distality was introduced by Hilbert (cf.\
Ellis \cite{El}) and studied by many in different contexts: see
Ellis \cite{El}, Furstenberg \cite{Fu} and Raja-Shah \cite{RS} and
the references cited therein.

Let $G$ be a locally compact (Hausdorff) group and let $e$ denote
the identity of $G$. Let $\Ga$ be a semigroup acting continuously on
$G$ by endomorphisms.  Then the $\Gamma$-action on $G$ is distal if
and only if $e\not\in\overline {\Gamma x}$ for all $x\in G\mi\{e\}$.
The group $G$ itself is said to be {\it distal} (resp.\ {\it
pointwise distal}) if the conjugacy action of $G$ on $G$ is distal
(resp.\ pointwise distal).

Abelian groups, discrete groups and compact groups are obviously distal.
It was shown by Rosenblatt in \cite{Ro} that the class of distal groups is closed
under compact extensions.  Nilpotent groups, connected
groups of polynomial growth are distal (cf.\ \cite{Ro}); recall that
a locally compact group $G$ with a left Haar measure $\lam _G$ is
said to be a group of {\it polynomial growth} if for each relatively
compact neighborhood $U$ of $e$ in $G$ there is a $k\in \N$ such
that $\{ {\frac{1}{n^k}}\lam _G(U^n) \mid n\geq 1 \}$ is bounded.  It
may also be noted that $p$-adic Lie groups of type $R$ and $p$-adic Lie
groups of polynomial growth are pointwise distal (cf.\ Raja
\cite{R1} and \cite{R2}): pointwise distal groups are called
noncontracting in Raja \cite{R1} and Rosenblatt \cite{Ro}.

The following types of groups
were studied by Losert in \cite{L1} and \cite{L2}.  A locally compact (resp.\ discrete) group $\Ga$ is called a 
generalized
$\ol{FC}$-group (resp.\ polycyclic) if $\Ga$ has a series $\Ga = \Ga_0 \supset \Ga_1
 \supset \cdots\supset \Ga_n = \{ e \}$ of closed normal subgroups such that
$\Ga_i/\Ga_{i+1}$ is a compactly generated group with relatively
compact conjugacy classes (resp.\ $\Ga_i/\Ga_{i+1}$ is cyclic) for
$i = 0, 1, . . . , n -1$. Note that polycyclic groups and compactly
generated groups of polynomial growth are generalized
$\ol{FC}$-groups. More detailed results on generalized
$\ol{FC}$-groups may be found in \cite{L2}.


Clearly, distal actions are pointwise distal and the converse holds
in case the action is on real Lie groups (see Abels [2] and Conze-Guivarsc'h [7]), but 
there are pointwise
distal actions which are not distal (see Jaworski-Raja \cite{JR} and
Rosenblatt \cite{Ro} for instance).  It can be easily seen that the
converse
can be extended in case $G$ and $\Ga$, which acts on $G$ by automorphisms, satisfy the following (see
Lemma \ref{ac}).  We say that a locally
compact group $G$ is {\it $\Ga$-Lie projective}
if every neighbourhood $U$ of the identity $e$ in $G$ contains a compact
normal $\Ga$-invariant subgroup $K$ such that $G/K$ is a
Lie group. In case $G$ is $\Ga$-{\it Lie projective} for $\Ga = \{ \ap ^n \mid n \in \Z \}$ for some $\ap\in\Aut(G)$,
we say that $G$ is $\ap$-{\it Lie projective}. Note that any $\Ga$-Lie projective group $G$ is Lie projective.
We show that the converse holds in case $\Ga$ is a generalized $\ol {FC}$-group whose action on
the group $G$ is pointwise distal (see Corollary \ref{lpc1}).

Proposition 1 of \cite{L2} shows that a generalized $\ol{FC}$-group
$G$ contains a compact normal subgroup $K$ such that $G/K$ is a Lie
group.  Additionally, if $G$ is pointwise distal, we show that it is Lie projective
(see Theorem \ref{FC}). This would
enable Lie theoretic considerations on generalized $\ol{FC}$-groups
that are pointwise distal. Moreover, this shows that any pointwise distal
generalized $\ol{FC}$ group is distal.

Distal actions of generalized $\ol{FC}$-groups were
considered in Jaworski-Raja \cite{JR}, Raja \cite{R3} and Shah \cite{Sh}. We quote some
related results: For the action of a polycyclic group $\Ga$ on a totally
disconnected group and for the action of a generalized $\ol
{FC}$-group on a compact metrizable group, distality and pointwise
distality are equivalent, (see Corollary 2.4 of \cite{JR} and Theorem
4.1 of \cite{R3}). We prove this equivalence for actions of generalized
$\ol{FC}$-groups on any locally compact group (Theorem \ref{distal-pd}).

We next look at the contraction group of automorphisms.  Let
$\ap\in \Aut(G)$, the group of bi-continuous automorphisms of $G$.
For a $\ap$-invariant compact subgroup $K$,
the $K$-{\it contraction group} of $\ap$ is defined by
$$
C_K(\ap ) = \{ g\in G\mid \ap ^n (g)K \to K \}.$$ The group $C_{\{
e\}}(\ap )$ is denoted by $C(\ap )$ and is called the {\it contraction
group} of $\ap$.  It is easy to see that $C(\ap ) = \{e\}$ if the
$\{\ap^n\}_{n\in\N}$-action on $G$ is distal.  It is an interesting
question to look at the converse: If the $\{\ap^n\}_{n\in\N}$-action
on $G$ is not distal, is it possible to find a $x \not = e$ such
that $\ap ^n (x) \to e$ as $n \to \infty$? Using exponential map of
Lie groups and results by Abels in \cite{Ab1} and \cite{Ab2}, the
converse can easily be seen to hold for (connected) Lie groups and,
it was proved in case of totally disconnected groups by Jaworski and Raja in
\cite{JR} using results of Baumgartner
and Willis from \cite{BW}. Recently the converse is also proved in
case of compact groups by Jaworski (cf.\ \cite{Jnew}),  We
show that the converse holds for all locally compact groups; namely,
we prove that  the $\{\ap^n\}_{n\in\N}$-action is distal if and only if
$C(\ap)=\{e\}$ (see Theorem \ref{ct}).

Main ingredient in the proof of the above result is
a decomposition of the type $C_K(\ap ) =C(\ap )K$ which is
proven under certain conditions (see Theorem \ref{clt-ctl}).
It is easy to see that $C(\ap )K\subset C_K( \ap)$.  The equality in
the above relation is a interesting question with applications to
probability (cf.\ \cite{HS}, \cite{DS1}) and this question was
answered affirmatively for real (resp.\ $p$-adic) Lie groups by Hazod
and Siebert in \cite{HS} (resp.\ by Dani-Shah in \cite{DS1}) and
for all totally disconnected groups by Baumgartner and Willis in
\cite{BW}.  Recently Jaworski has studied this
decomposition on compact groups in detail, (see \cite{Jnew} and related
references cited therein). In particular,
he illustrates by a counter example that such a decomposition
theorem may not hold, and it holds for an automorphism $\ap$ of a compact group 
$G$ if and only if it contains arbitrarily small closed normal $\ap$-invariant subgroups $N$
such that $G/N$ is finite dimensional (see Theorem 10.2 in
\cite{Jnew}).



We further explore properties of distal groups.  We first observe
that pointwise distal groups are unimodular (see Proposition \ref{pum});
recall that a locally compact
group $G$ is said to be {\it unimodular} if any left-invariant Haar measure is
also right-invariant.  The converse need not hold as any non-compact
semisimple connected Lie group is unimodular but it is not pointwise
distal. Observe that closed subgroups of pointwise distal groups are
also pointwise distal and hence, any pointwise distal group has the
property that all its closed subgroups are unimodular. We now
attempt to see how far one could progress in the converse direction.
Main result in Section 5 relates pointwise distality, contraction
groups and unimodularity of closed subgroups (see  Theorem
\ref{csu}).  It may be noted that such class of unimodular groups
arises in the classification of groups that admit recurrent random
walks (see \cite{GR} for any unexplained notions) and study of
contraction subgroups on such groups plays a crucial role in
\cite{GR}.  Here, we prove that the closed subgroups of $G$ are
unimodular if and only if for any inner automorphism $\ap$ and for
any $\ap$-invariant compact subgroup $K$ of $G$, $C_K(\ap)$ is
relatively compact (see Theorem \ref{csu}).

Let $G$ be a locally compact group and $P(G)$ denote the space of all
regular Borel probability measures on $G$ with weak* topology.  Let $\Aut(G)$ denote the group of
all bi-continuous automorphisms of $G$ and let $\Ga $ be a topological semigroup acting
continuously on $G$ (by automorphisms), i.e.\ the map $(\ga, g)\mapsto \ga(g)$
from $\Ga\times G$ to $G$ is continuous.  
Then this action extends to
a continuous action (of $\Ga$) on $P(G)$, which is given by
$\ap (\mu ) (E) = \mu (\ap^{-1} (E))$ for any $\ap \in \Ga$, any
$\mu \in P(G)$ and for any measurable subset $E$ of $G$.

The following proposition is useful in reducing the action on
general groups to that on metrizable groups.  In case $\Ga$ is countable
the result is proved in \cite{J}, \cite{Jnew}, (see the proof of Theorem 1 in \cite{J} 
and proof of Corollary 2.3 in \cite{Jnew}).  

\bp\label{am} Let $\Ga$ be a $\sigma$-compact locally compact group
acting on a $\sigma$-compact locally compact group $G$ by
automorphisms.  Then for each countable collection $\{ U_n \}$ of
neighborhoods of the identity $e$ in $G$ there exists a
$\Ga$-invariant compact normal
subgroup $L$ in $G$ such that $L\subset \cap U_n$ and $G/L$ is
metrizable with a countable basis for its open sets.\ep

The semidirect product, $\Gamma\ltimes G$ is also a $\sigma$-compact locally 
compact group and the assertion in the proposition easily follows from 
Theorem 8.7 of \cite{HR}. We will not give more details here. 

\end{section}

\begin{section}{Actions on Compact Groups}

In this section we discuss actions of semigroups on a compact group.
Let $(X, {\cal B}, m)$ be a probability space.   Recall that the action of a
semigroup of measure preserving transformations $\Ga$ on $X$ is said
to be {\it ergodic} if for any $\Ga$-invariant set $B\in \cal B$ (that is,
$\ap ^{-1}(B)=B$ for all $\ap \in \Ga$),  we have
$m(B)=0$ or $1$. For a compact group $K$ and $\Gamma\subset
\Aut(K)$, the $\Ga$-action on the homogeneous space $K/H$ for any
closed $\Ga$-invariant subgroup $H$, ergodicity is defined with
respect to the $K$-invariant probability measure on $K/H$.

Throughout this section, let $K$ denote a compact group and let
$\Gamma$ denote a topological semigroup acting (continuously) on $K$
by automorphisms.

The following is a generalization of Proposition 2.1 of \cite{R3} to
not necessarily metrizable groups for semigroup actions.
Note that for an action of a countable
group on a metrizable compact group, the result was already known
due to the work of Kitchens and Schmidt (see Theorem 2.3 in \cite{KS}
and Theorem 1.4 in \cite{S}), and in that case, the condition
of metrizability was removed by Jaworski (see Theorem 2.6 in \cite{Jnew}).

\bp\label{non-metric} There exists smallest closed $($resp.\
closed normal$)$ $\Gamma$-invariant subgroup $C$ $($resp.\ $C_1)$ in $K$
such that the $\Gamma$-action on $K/C$ $($resp.\ $K/C_1)$ is distal.
Moreover, the $\Gamma$-action on $C$ $($resp.\ $C_1)$ is ergodic and
$C\subset C_1$.

Moreover, if $K$ is metrizable, then $C$ is normal in $K$ and $C$ is also
the largest closed $\Ga$-invariant subgroup such that the $\Ga$-action on
$C$ is ergodic. \ep

\bo Without loss of any generality, we may assume that $\Ga\subset
\Aut(K)$.

\medskip
\noindent{\bf Step 1} Let $\K$ be the set of compact subgroups $L$ of $K$ such
that the $\Ga$-action on $K/L$ is distal. Here $\K$ is nonempty as
$K\in\K$. Let $C=\cap \{L\mid L\in\K\}$. Then $C$ is a compact
$\Ga$-invariant subgroup. Now we show that $\Ga$ acts distally on
$K/C$, i.e.\ $C\in\K$. Let $xC\in K/C$. Suppose $C$ belongs to the
closure of $\{\ga(x)C \mid \ga\in \Ga\}$ in $K/C$. Then for every
$L\in\K$, since $C\subset L$, we have that $L$ belongs to the
closure of $\{\ga (x)L\mid\ga\in \Ga\}$ in $K/L$ and since the
$\Ga$-action on $K/L$ is distal, $x\in L$, and hence $x\in C$.
Therefore, the $\Ga$-action on $K/C$ is distal.

\medskip
\noindent{\bf Step 2}
We show that the $\Ga$-action on $C$ is ergodic. If possible, suppose
it is not ergodic. By Theorem 2.1 of \cite{Be},
there exists an irreducible unitary representation $\pi$ of $C$ such that
$\pi_\Gamma=\{\pi\circ\gamma\mid\gamma \in \Gamma\}$ is finite up to
unitary equivalence.  Then there exist
$\gamma_1,\ldots,\gamma_n\in\Gamma$ such that for any
$\gamma\in\Gamma$, $\pi\circ\gamma$ is unitarily equivalent to
$\pi\circ\gamma_i$ for some $i$.  This implies that the
finite-dimensional unitary representation $\tilde \pi = \oplus \pi
\circ \gamma _i$ is $\Gamma$-invariant.  Let $C'= \Ker (\tilde \pi)$
and let $n$ be the dimension of $\tilde \pi$.  Observe that $C'\ne C$. Suppose
$\gamma_d(x)C'\to C'$ in $C/C'$ for some $x\in C$. Then $\tilde \pi
(\gamma_d(x))\to I_n$ where $I_n$ is the trivial operator on $\C
^n$.  Since $\tilde \pi$ is $\Gamma$-invariant, $\tilde
\pi\circ\gamma_d=u_d^{-1}\tilde \pi u_d$ for some $u_d\in U_n(\C)$,
hence
$$\tilde \pi(\gamma_d(x))=u_d^{-1}\tilde \pi(x)u_d \to I_n.$$
Since $U_n(\C)$ is compact, we get that $\tilde \pi(x)=I_n$, that is
$x\in \Ker (\tilde \pi)=C'$. This shows that the $\Ga$-action on
$C/C'$ is distal and hence it is distal on $G/C'$. This contradicts the minimality of $C$. 
This implies that the $\Ga$-action on $C$ is ergodic.

\medskip
\noindent{\bf Step 3}
Now we want to show the existence of the smallest closed normal
$\Ga$-invariant subgroup $C_1$ such that the $\Ga$-action on $K/C_1$ is
distal and the $\Ga$-action on $C_1$ is ergodic.  Take
$\Ga_1=\Ga.\Inn(K)$. This is a closed subsemigroup of $\Aut(K)$ as
$\Inn(K)$ is a normal compact subgroup of $\Aut(K)$.  Then from
above we get that there exists a smallest
$\Ga_1$-invariant subgroup $C_1$ of $K$ such that the $\Ga_1$-action on
$K/C_1$ is distal and the $\Ga_1$-action on $C_1$ is ergodic.  Here,
$C_1$ is normal in $K$ as $\Inn(K)\subset \Gamma_1$. Also, $\Ga$
acts distally on $K/C_1$.  Moreover, any $\Ga$-invariant normal
subgroup is also $\Ga_1$-invariant and hence $C_1$ is the smallest
closed normal $\Ga$-invariant group such that $\Ga$ acts
distally on $G/C_1$.  It is easy to verify that for any inner
automorphism $\ap$ defined by $x \in K$, $\pi \circ \ap$ is
unitarily equivalent to $\pi$, where equivalence is given by
$\pi(x)$.  This implies that the $\Ga$-action on $C_1$ is ergodic (see
also Lemma 2.4 of \cite{R3}). Also, the group $C$ as above is
contained in $C_1$.

\medskip
\noindent{\bf Step 4} Now suppose $K$ is metrizable.  Then there exits a dense $\Ga$-orbit
in $C_1$, and hence, in $C_1/C$.
But since the $\Ga$-action on $C_1/C$ is distal, we have that
$C_1=C$ and hence $C$ is normal. Let $C_2$ be a $\Ga$-invariant
subgroup of $K$ such that the $\Ga$-action on $C_2$ is ergodic. We need
to show that $C_2\subset C$. Then
the $\Ga$-action on $C_2C/C$ is also ergodic and hence it has a dense
$\Ga$-orbit. But the $\Ga$-action on $K/C$ is distal and hence
$C_2\subset C$. \hfill{$\Box$}\eo

\br \label{rem1} 1. As we can assume that $\Ga$ is a semigroup contained
in $\Aut(K)$, let $[\Ga]$ be the subgroup generated by $\Ga$ in
$\Aut(K)$, From Theorem 1 of \cite{El}, it is obvious that the
$\Ga$-action on $K$ is distal if and only if the $[\Ga]$-action on
$K$ is distal (see also the proof of Theorem 3.1 in \cite{RS}).  By
Theorem 2.1 of \cite{Be}, the same statement holds for ergodic
actions on a compact group. In this case we get in Proposition
\ref{non-metric} that the $\Ga$-action on $K/C$ (resp.\ on $C$) is
distal (resp.\ ergodic) iff the $[\Ga]$-action on $K/C$ (resp.\ on
$C$) is distal (resp.\ ergodic).

2. If $\Gamma$ (and hence $[\Gamma]$) is compactly generated then
metrizability of $K$ is not essential in the second assertion
because in this case, by Proposition \ref{am} we get that $K$ has
arbitrarily small compact normal $\Ga$-invariant subgroups $K_d$
such that $K/K_d$ is metrizable and we can argue as above for the
$\Ga$-action on each $K/K_d$ and get the desired assertion for $K$. \hfill{$\Box$}\er

\bc \label{cor1} Suppose $K$ is metrizable and suppose $K$ has a proper closed $\Ga$-invariant 
subgroup $L$. Then the $\Ga$-action on $K/L$ is distal if
and only if the $\Ga$-action is not ergodic on $H/C$, for any pair of compact
$\Ga$-invariant subgroups  $C$ and $H$ of $K$ such that $L\subset
C\subset H$ and $C\ne H$. \ec


\bo
Suppose there exist $\Gamma$-invariant subgroups $H$ and $C$ with
$L\subset C\subset H$ and $C\ne H$ such that the $\Ga$-action on
$H/C$ is ergodic. Then as $K$ is metrizable, so is $H/C$ and hence
there exists a dense $\Ga$-orbit in $H/C$, (see Theorem 5.6 of
\cite{Wa} for single transformation and the proof works for any
semigroup action also).
This in turn implies that the $\Ga$-action
on $H/C$ is not distal.  The canonical projection
$H/L\to H/C$ implies that $H/C$ is a factor of $H/L$.
Since any factor action of
a distal action is distal, the $\Ga$-action is not distal on $H/L$,
hence on $K/L$ (cf.\ Corollary 6.10 of \cite{BJM}).  

We now prove the converse. Suppose
the $\Ga$-action on $K/L$ is not distal. By Proposition
\ref{non-metric}, there exists a closed normal $\Ga$-invariant
subgroup $H'$ of $K$ such that the $\Ga$-action on $H'$ is ergodic
and it is distal on $K/H'$.  Let $H=H'L$ and $C=L$. Then the
$\Ga$-action on $H/C\simeq H'/H'\cap L$ is ergodic. Here $H\ne C$,
otherwise, $H'\subset L$ and since the $\Ga$-action on $K/H'$ is
distal, the $\Ga$-action on $K/L$ is also distal, a contradiction to
our assumption. \hfill{$\Box$}\eo

Recall that $(\Gamma,K)$ satisfies DCC ({\it descending chain
condition}) if for each sequence $\{ K_n \}_{n\in\N}$ of compact
$\Ga$-invariant subgroups of $K$ such that $K_n\supset K_{n+1}$,
$n\in\N$, there exists $n_0\in \N$ such that $K_n=K_{n_0}$ for all
$n\geq n_0$.

The following may be proved along the lines of Theorem
3.15 of \cite{KS} but here we give a different proof using
Proposition \ref{non-metric}. Note that it was already shown for a countable group
$\Ga$, see Proposition 3.5 in \cite{S}) when $K$ is metrizable, a
condition which was removed by Jaworski (cf.\ \cite{Jnew}, Theorem 2.6).

\bl\label{DCC-Lie} Suppose $K$ is metrizable and $(\Gamma, K)$ satisfies DCC.
Then there exists a compact normal $\Gamma$-invariant subgroup $C$
of $K$ such that $K/C$ is a real Lie group, the $\Gamma$-action on
$K/C$ is distal and the $\Gamma$-action on $C$ is ergodic. In particular,
if the $\Ga$-action on $K$ is distal, then $K$ is a Lie group. \el

\bo From Proposition \ref{non-metric}, there exists a compact normal
$\Gamma$-invariant subgroup $C$ such that the $\Gamma$-action on
$K/C$ is distal and the $\Gamma$-action on $C$ is ergodic.  As $K$ is metrizable,
$C$ is the largest subgroup such that the $\Ga$-action on $C$ is ergodic. It is
easy to verify that $(\Ga, K/C)$ also satisfies DCC.  Hence it is
enough to prove that if the $\Gamma$-action on $K$ is distal and
$(\Ga, K)$ satisfies DCC, then $K$ is a Lie group.

Since $K$ is a metric group and the $\Ga$-action on $K$ is distal, it is not
ergodic on any nontrivial subgroup of $K$ (see Proposition \ref{non-metric} or
Corollary \ref{cor1}).

Suppose $K$ is not a real Lie group.
As in Step 2 of Proposition \ref{non-metric}, since the $\Ga$-action on
$K$ is not ergodic, there exists a $\Ga$-invariant proper closed
non-trivial normal subgroup $K_1$ of $K$ such that $K/K_1$ is a Lie group.
Proceeding this way we obtain a strictly decreasing sequence
$\{ K_n \}$ of $\Ga$-invariant closed (normal) subgroups of $K$ such
that $K_0=K$.  This is a contradiction to the hypothesis that
$(\Gamma, K)$ satisfies DCC. Hence $K$ is a real Lie group. \hfill{$\Box$}\eo

\bl\label{ind} Suppose $\Ga$ is a subgroup of $\Aut(K)$ and there
exists a $\Ga$-invariant compact normal subgroup $C$ in $K$ such
that $K/C$ is a Lie group. Suppose also that $\gamma \in\Aut(K)$ is
such that $\gamma $ normalizes $\Ga$ and $\{\gamma^n\}_{n\in\Z}$
acts distally on $K$. Then there exists a compact normal
$\Ga$-invariant subgroup $C'$ contained in $C$ such that $\gamma(C')=C'$ and
$K/C'$ is a Lie group. In particular $C'$ is invariant under the
group generated by $\gamma$ and $\Gamma$. Moreover, there exist
compact normal $\ga$-invariant subgroups $C'_d$ such that
each $K/C'_d$ is a Lie group and $\cap_d C'_d=\{e\}$; that is, 
$K$ is $\ga$-Lie projective.   \el

\bo Let $C'=\cap_{n\in\Z}\gamma^n(C)$. It is clearly
$\gamma$-invariant and normal in $K$. Moreover, for a
$\gamma'\in\Ga$, and $n\in\Z$, let
$\gamma_n=\gamma^{-n}\ga'\gamma^n$ then $\ga_n\in \Ga$ and
$\ga'(\gamma^n(C))=\gamma^n(\gamma_n(C))=\gamma^n(C)$. Hence
$\ga'(C')=C'$, i.e.\ $C'$ is $\Ga$-invariant.

Let $G=K/C$. Then $G$ is a Lie group and we define a map
$\vp:K/C'\to G^\Z$ as follows: $\vp(gC')=(g_n)_{n\in\Z}$, where
$g_n=\ga^n(g)C\in G$. It is easy to see that $\vp$ is a continuous
bijective homomorphism from $K/C'$ onto $\vp(K/C')\subset G^\Z$, and
$\vp\circ\gamma=\alpha\circ\vp$, where $\alpha$ denotes the shift
map on $G^\Z$. Here, $(\{\alpha^n\}_{n\in\Z},G^\Z)$ satisfies DCC
(cf.\ \cite{KS}). Since $\{\ga^n\}_{n\in\Z}$ acts distally on
$K$, it also acts distally on $K/C'$ (cf.\ \cite{RS}, Theorem 3.1).
Now by Lemma \ref{DCC-Lie}, $K/C'$ is a Lie group.

As $K$ is Lie projective, there exist compact normal subgroups $C_d$ such
that $K/C_d$ is a Lie group and $\cap_d C_d=\{e\}$. Now from the first assertion
for $\Ga$ as the trivial subgroup of $\Aut(K)$, there exist $\ga$-invariant
compact normal subgroups $C'_d$ contained in $C_d$ such that $K/C'_d$ is a Lie group.
Clearly, $\cap_d C'_d=\{e\}$.
\hfill{$\Box$}\eo

\bp\label{pd-FC} Suppose $\Ga$ is a generalized $\ol{FC}$-group  which acts on $K$
by automorphisms. Suppose the $\Ga$-action on $K$ is pointwise distal. Then
$K$ is $\Ga$-Lie projective. \ep

\bo Let $K^0$ and $\Ga ^0$ be the connected component of the identity in
$K$ and $\Ga$ respectively. Observe that $\Gamma$ also acts on $\Inn(K)$ by automorphisms. 
We form the semidirect product
$\Ga'=\Ga\ltimes \Inn(K)$; which is algebraically the product of the groups $\Ga$ and $\Inn(K)$ and the
topology on $\Ga'$ is the product topology of $\Ga$ and the compact topology of $\Inn(K)$. Here,
$\Gamma'$ is also a generalized $\ol{FC}$-group. Moreover, it can be verified that $\Ga'$ acts on $K$ and the
action is pointwise distal.  Replacing $\Ga$ by $\Ga'$ if necessary, we
may also assume that $\Inn(K)$ is a compact subgroup $\Ga$. Let $\pi:\Gamma\to\Aut(K)$ be the natural continuous
homomorphism given by the action. Here $\ker\pi=\{\gamma\in\Gamma\mid\gamma(k)=k\hbox{ for all }k\in K\}$ is a closed
subgroup of $\Ga$ and since $\Ga$ is a generalized $\ol{FC}$-group, so is $\Gamma/\ker\pi$;
the latter group is algebraically isomorphic to $\pi(\Ga)$. Moreover, any
$\pi(\Ga)$-invariant group is $\Ga$-invariant. Replacing $\Ga$ by $\pi(\Ga)$ if necessary,
we may assume that $\Gamma$ is a subgroup of $\Aut(K)$ and $\Ga$ carries a topology finer than
the subspace topology with respect to which it is a generalized $\ol{FC}$-group and
$\Inn(K)$ is a compact subgroup $\Ga$.  It is known that for a compact group $K$,
the connected component of $\Aut(K)$ is $\Inn(K^0)$, the group of inner automorphisms
by elements of $K^0$ (cf.\ \cite{I}, Theorem $1'$). As $\Inn(K)\subset \Ga$, we get that $\Ga^0=\Inn(K^0)$. In
particular, $\Ga^0$ is compact. So, there exists a compact normal
subgroup $L$ in $\Ga$ such that $\Ga/L$ is discrete and has a
(normal) polycyclic subgroup of finite index (see Proposition 2.8 of
\cite{JR}). In particular, $\Ga/L$ is finitely generated.  Let
$F=\{\gamma_i\in\Ga\mid i=1,\ldots,n\}$ be such that the following
hold: there exist
$L=\Gamma_0\subset\ldots\subset\Gamma_n\subset \Ga$, each $\Gamma_i$ is a
normal subgroup of $\Gamma$, $\Ga/\Ga_n$ is finite and each $\Gamma_i$ is generated by
$\gamma_i$ and $\Gamma_{i-1}$, $1\leq i\leq n$.

Here, $L\ltimes K$ is compact and hence Lie projective and there
exist compact normal subgroups $H_d$ in $L\ltimes K$ such that
$\cap_d H_d=\{e\}$ and each $(L\ltimes K)/H_d$ is a Lie group. Let
$C_d=K\cap H_d$. Then each $C_d$ is a $L$-invariant normal subgroup
in $K$ such that $K/C_d$ is a Lie group and $\cap_d C_d=\{e\}$. Applying
Lemma \ref{ind} successively for $\Ga _ i$ and $\gamma_{i+1}$,
$0\leq i\leq n-1$, we get that there exist compact normal $\Ga_n$-invariant
subgroups $M_d\subset C_d$ such that $K/M_d$ is a Lie group. Since
$\Ga/\Ga_n$ is finite, $K_d=\cap_{\ga\in\Ga}\gamma(M_d)$ is a finite intersection
and hence $K/K_d$ is a Lie group.
Clearly, $K_d$ is $\Ga$-invariant and $\cap_d K_d=\{e\}$. \hfill{$\Box$}\eo


\end{section}

\begin{section}{Distal and Pointwise Distal Groups}

In this section we compare distality and pointwise distality of
certain actions. We know from Rosenblatt \cite{Ro} that distality,
pointwise distality and polynomial growth are all equivalent
properties for any almost connected group. This is not true in general as
there are abelian extensions of compact groups that are pointwise
distal but not distal (see Example 2.5 of \cite{JR}). There are also
examples of $\Z$-extensions of compact groups which are not
poinwise distal. Our first result shows a relation between
polynomial growth and distality for groups which need not be almost connected.


\bt\label{pgrowth} Let $G$ be a compactly generated locally compact
group of polynomial growth. Then $G$ has a compact normal subgroup
$C$ such that $G/C$ is distal and the conjugacy action of $G$ on $C$
is ergodic.  \et

\bo We first assume that $G$ is Lie group such that $G^0$ has no
nontrivial compact normal subgroup. Since $G/G^0$ is discrete, it is
enough if we show that the conjugacy action of $G$ on $G^0$ is distal. For each
$g\in G$, let $\alpha_g$ denote the automorphism of $G^0$ defined by
the conjugation action of $g$ restricted to $G^0$ and let $\du\alpha_g$
denote the corresponding Lie algebra automorphism of the Lie algebra
${\cal G}$ of $G^0$. Since $G^0$ has no nontrivial compact normal subgroup,
Theorem 1 of \cite{L1} implies that the eigenvalue of
$\du\alpha_g$ are of absolute value 1.  Hence by Theorem 1 of
\cite{Ab1} and Theorem 1.1 of \cite{Ab2}, the conjugacy action of
$G$ on $G^0$ is distal.

Suppose $G$ is not a Lie group. Then $G$ has a maximal compact
normal subgroup $K$ such that $G/K$ is a Lie group (see \cite{L1}).
Since $(G/K)^0$ has no nontrivial compact normal subgroup, from
above, $G/K$ is distal. By Proposition \ref{non-metric} and Remark
\ref{rem1} (2), we have that $K$ has a maximal compact $G$-invariant
(normal) subgroup $C$ such that the conjugacy action of $G$ on $K/C$, and
hence, on $G/C$ is distal and the conjugacy action of $G$ on $C$ is
ergodic. \hfill{$\Box$}\eo

We know that distal groups are pointwise distal and we also know
that the converse is not true (see for instance, Example 2.5 of
\cite{JR}). Observe that if $G$ is a (not necessarily connected) Lie group,
and if the action of a group $\Ga$ on $G$ is pointwise distal then it is distal,
(this follows from Theorem 1.1 of \cite{Ab2} and \cite{CG}).
Hence the following also holds.

\bl\label{ac}

Let $G$ be a locally compact group and let $\Ga$ be a group which acts on $G$
by automorphisms such that $G$ is $\Ga$-Lie projective.  Then the $\Ga$-action on $G$ is distal
if and only if it is pointwise distal. In particular, any Lie-projective pointwise distal group is
distal. A pointwise distal group which is a closed subgroup of an almost connected group is distal.
\el

We now compare distality and pointwise distality of a particular
class of locally compact groups.

\bt\label{FC} Let $G$ be a generalized $\ol{FC}$-group.  Suppose
$G$ is pointwise-distal. Then $G$ is Lie projective.  Moreover, $G$
is distal. \et

\bo Due to Lemma \ref{ac},
it is sufficient to prove that any pointwise distal
generalized $\ol{FC}$-group is Lie projective.

Suppose $G$ is a generalized $\ol{FC}$-group.  Then there exists a
maximal compact normal subgroup $K$ in $G$ such that $G/K$ is a Lie
group (cf.\ \cite{L2}).  Since the conjugacy action of $G$ on $K$ is
pointwise distal, by Proposition \ref{pd-FC} there exist closed
normal subgroups $K_d$ of $K$, which are invariant under the conjugacy action
of $G$, such that each $K/K_d$ is a Lie group and $\cap_d K_d=\{e\}$.
That is, each $K_d$ is normal in $G$ and since
$G/K$ is isomorphic to $(G/K_d)/(K/K_d)$, where $G/K$ and $K/K_d$ are Lie groups, we
get that $G/K_d$ is a Lie group for every $d$. This shows that $G$
is Lie projective. \hfill{$\Box$}\eo

Since any compactly generated locally compact group of polynomial
growth is a generalized $\ol{FC}$-group (see \cite{L2}) we have the
following:

\bc\label{lie-proj} Let $G$ be a compactly generated locally compact
pointwise-distal group of polynomial growth. Then $G$ is Lie
projective. Moreover, $G$ is distal. \ec

\br \label{rem2} From Theorem \ref{pgrowth} and Corollary
\ref{lie-proj}, it follows that any compactly generated locally
compact group of polynomial growth has a compact normal subgroup $C$
such that $G/C$ is distal and Lie projective and the conjugacy
action of $G$ on $C$ is ergodic. \hfill{$\Box$}\er

We next compare distality and pointwise distality of an action
of a generalized $\ol{FC}$-group on a locally compact group. Note that
it is proved in case of metrizable compact groups in
\cite{Ra} by a different method.

\bt\label{distal-pd} Let $G$ be a locally compact group and let
$\Gamma$ be a generalized $\ol{FC}$-group acting on $G$ by
automorphisms. Then the $\Gamma$-action on $G$ is distal if and only
if the $\Gamma$-action on $G$ is pointwise distal. \et

\bo One way implication ``only if'' is obvious. Now suppose the
$\Gamma$-action on $G$ is pointwise distal. Let $G^0$ be the
connected component of the identity in $G$.  Then there is a maximal
compact normal subgroup $K$ of $G^0$ such that $G^0/K$ is a Lie
group. Since $K$ is maximal, $K$ is characteristic in $G$ and hence
$\Ga$-invariant. Since the $\Ga$-action on $K$ is pointwise distal.
By Proposition \ref{pd-FC}, $K$ has compact normal $\Ga$-invariant
subgroups $K_d$ such that $\cap K_d=\{e\}$ and $K/K_d$ is a Lie
group. Since the $G^0$-action on $K$ is by inner automorphisms of
$K$, each $K_d$ is normal in $G^0$ and hence $G^0/K_d$ is a Lie
group. By Theorem 3.1 of \cite{RS}, the $\Ga$-action on each
$G^0/K_d$ is pointwise distal and hence, it is distal by Theorem 1 of \cite{Ab1} and
Theorem 1.1 of \cite{Ab2}. Since this is true for each $d$ and $\cap_d
K_d=\{e\}$, we get that the $\Ga$-action on $G^0$ is distal.

Now it is enough to prove that the $\Ga$-action on $G/G^0$ is distal. We know that the
$\Gamma$-action on $G/G^0$ is pointwise distal by Theorem 3.3 of
\cite{RS}. Hence we may assume that $G$ is a totally disconnected
group. Since the connected component $\Ga^0$ of $\Ga$ acts
trivially on $G$ and $\Ga/\Ga^0$ is also a generalized
$\ol{FC}$-group and its action on $G$ is pointwise distal, we may assume
that $\Ga$ is totally disconnected.  By \cite{L2}, $\Ga$ contains a
compact normal subgroup $L$ such that $\Ga /L$ is a Lie group. Since
$\Ga$, and hence, $\Ga/L$ is also totally disconnected, $\Ga /L$ is discrete.
As $\Ga/L$ is a discrete generalized $\ol{FC}$-group, $\Ga /L$ has a polycyclic
subgroup of finite index (cf.\ \cite{L2}), Hence by Lemma 2.3 of
\cite{JR} (which is also valid for any non-metric group by \cite{J},
see also \cite{JR}), we get that the $\Ga$-action of $G$ is distal. \hfill{$\Box$}\eo








We next extend Proposition \ref{pd-FC} to actions on Lie projective groups.

\bc\label{lpc1}
Let $\Ga$ be a generalized $\ol{FC}$-group which acts on a locally compact group
$G$ by automorphisms. Suppose the $\Ga$-action on $G$ is pointwise distal and 
$G$ is Lie projective.  Then $G$ is $\Ga$-Lie projective.

\ec

\bo By Theorem \ref{distal-pd}, the $\Ga$-action on $G$ is distal. Now by Theorem 3.3 of \cite{RS},
the $\Ga$-action on $G/G^0$ is distal, and hence, it is equicontinuous (cf.\  \cite{Sh}, Proposition 1.5). I.e.\ every
neighbourhoold of the identity in $G/G^0$ contains a compact open $\Ga$-invariant subgroup.
In particular, there exists an almost connected open $\Ga$-invariant subgroup $H$ in $G$.
Let $K$ be the largest compact normal subgroup of $H$.
Then $K$ is characteristic in $H$ and hence $\Ga$-invariant. By Proposition \ref{pd-FC}, $K$ is $\Ga$-Lie projective.

Let $U$ be any neighbourhood of the identity $e$ in $G$. Then $U\cap K$ is open in $K$ and as $K$ is 
$\Ga$-Lie projective, there exists a compact $\Ga$-invariant subgroup $L\subset U\cap K$ such that $L$ is 
normal in $K$ and $K/L$ is a Lie subgroup. Let $V$ be a neighbourhood of the identity $e$ in $G$ such that 
$L\subset V\subset U$ and $V\cap K$ does not contain any compact subgroup larger than $L$. As $G$ is Lie 
projective, there exist a compact normal subgroup $M$ in $G$ such that  $M\subset V\cap H$ and $G/M$ is a Lie group. 
Here, $M\subset K$ as $K$ is the largest compact normal subgroup of $H$. Hence $M\subset V\cap K$. 
Now from above, we get that $M\subset L$. Since $L$ is $\Ga$-invaraint, we have that $\ga(M)\subset L$ for all 
$\ga\in \Ga$. Moreover, since $M$ is normal in $G$, so is $\ga(M)$ for all $\ga\in\Ga$. Let $M_0$ be the closed 
subgroup generated by $\ga(M)$, $\ga\in\Ga$. Then $M\subset M_0\subset L\subset U$ and $M_0$ is a compact normal 
$\Ga$-invariant subgroup of $G$. Moreover, $G/M_0$ is a Lie group as it is isomorphic to $(G/M)/(M_0/M)$. Since this holds for any such $U$ as above, $G$ is $\Ga$-Lie projective. \hfill{$\Box$} \eo
\end{section}

\begin{section}{Distality and Contraction Groups}

In this section we get a necessary and sufficient condition for
distality of the action of $\{\ap^n\}_{n\in\N}$ on a locally compact group for
an automorphism $\ap$ in terms of its contraction group, (see Theorem \ref{ct}).
We also get a decomposition theorem for contraction groups, (see Theorem \ref{clt-ctl}).

Recall that for an $\ap$-invariant subgroup $K$ of $G$,
the  $K$-contraction group of $\ap$ is defined as $C_K(\ap ) = \{ x \in G \mid
\ap ^n (x)K \to K ~~{\rm as}~~ n \to \infty \}$ and we denote
$C_{\{e \}}(\ap )$ by $C(\ap)$ which is known as the contraction group of $\ap$.

It is evident that for any automorphism $\ap$, if the
$\{\ap ^n\}_{n\in\N}$-action is distal, then $C(\ap)$ is trivial. But the converse is not
known except in case of compact groups which has been recently proved (cf.\ Theorem 13.2, \cite{Jnew}). We
prove the converse for all locally compact groups in this section, (see Theorem \ref{ct}).
We recall that an automorphism $\ap$ is distal on $G$ (equivalently, the $\ap$-action on $G$ is distal) if the 
$\{\ap^n\}_{n \in \Z}$-action on $G$ is distal.
Note that for a connected Lie group $G$, using Proposition 3.2.6 of \cite{HSb}, it follows from Abels (\cite{Ab1}, \cite{Ab2}),
that $\ap$ is distal on $G$ if and only if $C(\ap^{\pm 1})=\{e\}$;
the same holds also for totally disconnected groups by Jaworski and Raja (cf.\ \cite{JR}). The following result strengthens the
result by Abels and the one by Jaworski and Raja.

\bt \label{ct} Let $G$ be a locally compact group and
$\ap\in\Aut(G)$. Then the $\{ \ap^n \}_{n\in \N}$-action on $G$ is
distal if and only if $C(\ap )=\{ e \}$. In particular, $\ap$ is
distal if and only if $C(\ap ^{\pm 1}) = \{ e \}$.\et

To prove Theorem \ref{ct}, we
need to prove a decomposition theorem for contraction groups when
$G^0$ is $\ap$-Lie projective.

\bt\label{clt-ctl}
Let $G$ be a locally compact group and $ \ap\in\Aut(G)$.
Suppose $G^0$ is $\ap$-Lie projective.  Then $C_L(\ap)=C(\ap)L$
for any $\ap$-invariant compact subgroup $L$ of $G$.

In particular, $C_L(\ap)=C(\ap)L$ if  either $G^0$ is a Lie group or
if the restriction of $\alpha$ to the largest compact normal
subgroup of $G^0$ is distal. \et

The above decomposition theorem generalizes results proven for (connected)
Lie groups (cf.\ \cite{HS}) and for totally disconnected groups
(cf.\ \cite{BW}, \cite{J}). However, it is not true for all compact groups if the
automorphism is not distal; see \cite{Jnew} for a counter example. Note that the last condition above is
equivalent to the condition that the $\{\ap^n\}_{n\in\N}$-action on the maximal compact subgroup of $G^0$ is distal.

Towards the proof of the decomposition theorem, we prove some preliminary
results.

\bp \label{conn-ct} Let $G$ be a connected locally compact group,
$\ap\in\Aut(G)$ and let $K$ be the maximal compact normal subgroup
of $G$.  Suppose the $\{\ap^n\}_{n\in\N}$-action on $K$ is distal.
Then $C(\ap)$ is closed and a simply connected nilpotent group and for any compact $\ap$-invariant subgroup
$L$ of $G$, $C_L(\ap)=L\ltimes C(\ap)$ and $C_K(\ap)=C(\ap)\times K$. \ep

Note that from \cite{Si}, any closed connected contraction subgroup $C(\ap)$ is a simply connected
nilpotent Lie group. Before proving the proposition, we need the following lemma, (see Proposition 10 in \cite{J} for a similar result).

\bl\label{pl} Let $G$ be a locally compact group and let $\alpha\in\Aut(G)$. Suppose for an index set $I$,
there is a family $\{ K_i\}_{i\in I}$ of compact $\ap$-invariant subgroups of $G$ such that it is closed under finite intersection and
$\cap_i K_i = \{ e \}$
The the following hold:
\begin{enumerate}
\item[$(1)$] $\cap_i C_{K_i}(\ap)=C(\ap)$.

\item[$(2)$] If $L$ is a compact $\ap$-invariant subgroup of $G$ such that $C_L(\ap) \subset C_{K_i}(\ap)L$, then
$C_L(\ap ) = C(\ap )L$, if one of the following conditions are satisfied:
\begin{enumerate}
\item[$(a)$] The $\ap$-action on $L$ is distal.
\item[$(b)$] $C_{K_i}(\ap)\cap L\subset K_i$ for all $i\in I$.
\end{enumerate}
\end{enumerate}
\el

Note that conditions (a) and (b) in (2) are equivalent. For each $i$, as $L\cap K_i$ is $\ap$-invariant, if
the $\ap$-action on $L$ is distal, so is the corresponding action of $\ap$ on $L/L\cap K_i$ (cf.\  \cite{RS},
Theorem 3.1), and hence $C_{K_i}(\ap)\cap L=\C_{K_i\cap L}(\ap)\cap L\subset K_i$.  Conversely, suppose (b) holds.
As $C(\ap)\cap L\subset C_{K_i}(\ap)\cap L\subset K_i$ and $\cap_i K_i=\{e\}$, we have $C(\ap)\cap L=\{e\}$. Therefore, by
Theorem 13.2 of \cite{Jnew}, the $\ap$-action on $L$ is distal.

\bo $\!\!\!\!\!${\bf ~of Lemma \ref{pl}~} Let $\{K_i\}_{i\in I}$ be as in the hypothesis.
We first show that $\cap_i C_{K_i}(\ap)=C(\ap)$.  One way inclusion is obvious.
We may also assume that every $K_i$ is nontrivial. Let $K=K_i$ for a fixed $i\in I$. Then we can replace
$K_i$ by $K_i\cap K$ for every $i\in I$ and assume that each $K_i$ is
contained in a compact group $K$. Let $a\in \cap_iC_{K_i}(\ap)$. Given any neighbourhood $U$ of $e$,
there exists $j\in I$ such that $K_j\subset U$. If every $B_i=K_i\cap (G\setminus U)\ne\emptyset$. We know that
$\{B_i\}_{i\in I}$ is a collection of nonempty compact sets in $K$ which has finite intersection property as $\{K_i\}_{i\in I}$ is
closed under finite intersection.
Hence $B=\cap_iB_i$ is nonempty. As $B\subset \cap_i K_i=\{e\}$ and $B\subset G\setminus U$, we
arrive at a contradiction. As $a\in C_{K_j}(\ap)$ and $K_j\subset U$, we get that $\ap^n(a)\in U$ for all large $n$.
This implies that $a\in C(\ap)$. This proves (1)

As (a) and (b) are equivalent, we prove the assertion in (2) assuming (b).
Let $g\in C_L(\ap)$. It is enough to show that $gL\cap C(\ap)$ is
nonempty. From the hypothesis $g=x_il_i$ for some $x_i\in
C_{K_i}(\ap)$ and $l_i\in L$. Let $C_i=gL\cap C_{K_i}(\ap)$.  It is nonempty as it contains $x_i$.
Note that $x_i(L\cap K_i)\subset C_i$ and
$C_i=gL\cap C_{K_i}(\ap)=x_iL\cap C_{K_i}(\ap)=x_i(L\cap C_{K_i}(\ap))\subset L\cap K_i$ due to (b).
Hence $C_i=x_i(L\cap K_i)$ is compact for each $i$. Moreover, $C_i\cap C_j\supset gL\cap C_{K_i\cap K_j}(\ap)$. As
$\{K_i\}$ is closed under finite intersection, $\{C_i =x_i(K_i\cap L)\}_{i\in I}$ has finite intersection property;
it is a family of compact subsets of
$gL$. As $gL$ is compact, $C=\cap C_i$ is nonempty. That is,
$C=(\cap_i C_{K_i}(\ap))\cap gL=C(\ap)\cap gL$ is nonempty. \hfill{$\Box$}\eo

\bo  $\!\!\!\!\!${\bf ~of Proposition \ref{conn-ct}} Let $G$, $K$ and $\ap$ be as in the hypothesis. 
We know that $G/K$ is a connected Lie group without any compact central
subgroup of positive dimension.  Let $\ap'$ be the automorphism on
$G/K$ induced by $\ap$ which is defined as $\ap'(xK)=\ap(x)K$ for
all $x\in G$. We know that $C(\ap')$ is a connected nilpotent subgroup (cf.\ \cite{Si}),
and hence so is $\overline{C(\ap')}$. We first show that $C(\ap')$ is closed.
Let $D=\{d_n\}_{n\in\N}$ be a countable subgroup of $C(\ap')$ which is dense in
$\overline{C(\ap')}$. Let $\mu=\sum_{n\in\N} 2^{-n}\delta_{d_n}$. It is a probability measure
on $G/K$. Then it is easy to see that $(\ap')^n(\mu)\to\delta_e$. Therefore, by
Theorem 1.1 of \cite{DS2}, $\supp\mu\subset C(\ap')$. But
$\supp\mu=\overline{D}=\overline{C(\ap')}$. Hence
$C(\ap')$ is closed in $G/K$ and it is a simply
connected nilpotent group (cf.\ \cite{Si}). Therefore, $C_K(\ap)$ is
closed in $G$. Since the $\{\ap^n\}_{n\in\N}$-action, and hence,
the $\{\ap^n\}_{n\in\Z}$-action on $K$ is distal, $C(\ap)\cap K$ is
trivial.

We first assume that $G$ is a Lie group. Then $C_K(\ap ) = KC(\ap)$, more generally, 
$C_L(\ap)=LC(\ap)$ for any compact $\ap$-invariant subgroup $L$, (cf.\ \cite{HS}, Theorem 2.4). 
From above we have that $C_K(\ap)$ is closed. Moreover,
since the $\{\ap^n\mid n\in\N\}$-action on $K$ is distal, $C(\ap)\cap K=\{e\}$. Now
by Corollary 2.7 of \cite{HS}, $C(\ap)$ is closed. Moreover, $C(\ap)$ is a simply connected nilpotent group, 
it does not contain any non-trivial compact subgroup. 


Suppose $G$ is a connected (not necessarily Lie) group. Since the
action of $\{\ap^n\}_{n\in\Z}$ on $K$ is distal, by Lemma \ref{ind}
there exist $\ap$-invariant compact normal subgroups $K_i$ of $K$
such that $\cap_i K_i=\{e\}$ and $K/K_i$ are Lie groups (see also Theorem 2.6 and Corollary 2.7 in \cite{Jnew}).
Without loss of any generality, we may assume that $\{K_i\}$ is
closed under finite intersection. Since $G$ is connected and $K$ is a
compact normal subgroup of $G$, the action of $G$ on $K$ is by inner
automorphisms of $K$. This implies that $K_i$'s are normal in $G$.
Since $G/K$ and $K/K_i$ are real Lie groups, $G/K_i$'s are connected
real Lie groups. From above, as $\ap$-action on $K/K_i$ is also distal (cf.\ \cite{RS}, Theorem 3.1), $C_{K_i}(\ap)$ is
closed for each $i$. By Lemma 4.4 (1), $C(\ap)$ is also closed. Since the $\{\ap^n\}_{n\in\N}$-action on
$K$ is distal, by Lemma \ref{pl} above, $C_K(\ap)=KC(\ap)$. Let $L$ be any $\ap$-invariant subgroup of $G$. 
As $G/K$ is a Lie group
from Theorem 2.4 of \cite{HS}, we have that $C_L(\ap)\subset C_{LK}(\ap)=C_K(\ap)L=C(\ap)KL$.
Hence as $K$ is normal,
$C_L(\ap)=C(\ap)LK\cap C_L(\ap)=C(\ap)L(K\cap C_{K\cap L}(\ap))=C(\ap)L(K\cap L)=C(\ap)L$ (cf.\ \cite{RS}, Corollary 3.2).

As noted above $C(\ap)$ is closed, and hence,  
so is $C(\ap)L=C_L(\ap)$. As each $C_{K_i}(\ap)/{K_i}=C(\ap)K_i/K_i$ is a closed and simply connected 
(nilpotent) group of $G/K_i$ for $\{K_i\}$ as above,
$C(\ap)$ is also connected and, as it is closed, it is a simply connected and nilpotent group (cf.\ \cite{Si}). In particular, 
$C(\ap)\cap L= \{ e \}$. Therefore, $C_L(\ap)=L\ltimes C(\ap)$. Observe that
as $K$ is normal in $G$, we get that $C_K(\ap)=K\times C(\ap)$, a direct product. \hfill{$\Box$}\eo

\bl\label{lc} Let $G$ be a locally compact group and let
$\ap\in\Aut(G)$. Let $M$ be a compact $\alpha$-invariant subgroup of
$G^0$ such that it is normal in $G^0$ and $G^0/M$ is a Lie group.
Then for every neighbourhood $V$ of $e$ in $G$, $VG^0$ contains a
compact subgroup $C$ normalized by $G^0$ such that $CG^0$ is open,
$C\cap G^0=M$, $CG^0/M=C/M\times G^0/M$ and if $c\in C$ is such that
$\ap(c)\in CG^0$, then $\ap(c)\in C$. \el

\bo
Since $G^0/M$ is a Lie group, there exists
a relatively compact neighbourhood $W$ of $e$ such that $W\subset V$ and $WM\cap G^0$ does
not contain any compact subgroup larger than $M$. Let $U\subset W$ be a neighbourhood of $e$ in
$G$ such that $U\ap(U)\subset W$. Moreover, Since $U$ is open, so is $UG^0$ and it contains an open
Lie projective subgroup, say $H$. Since $U\cap H$ is open in $H$, it contains a compact
subgroup $C'$ which is normal in $H$ and $H/C'$ is a Lie group. That is, $C'$ is normalised by $G^0$ and
$C'G^0$ is open in $H$, and hence, in $G$. Let $C=C'M$. Observe that $C$ is a compact group normalised by 
$G^0$ and $CG^0=C'G^0$ is open. As $C\cap G^0\subset UM\cap G^0\subset WM\cap G^0$ and as $M\subset C$,
we have $C\cap G^0=M$. Therefore, $C$ normalises $M$, i.e.\ $M$ is normal in $CG^0$ and
$CG^0/M=C/M\times G^0/M$.

Let $\pi: CG^0\to CG^0/M$ be the
natural projection. Then  from above, $\pi(CG^0)$ is a direct product of $\pi(C)$ and $\pi(G^0)$, where
$\pi(C)$ is a compact totally disconnected group and $\pi(G^0)$ is a connected Lie group.

Let $c\in C$ be such that $\ap(c)\in CG^0$. As $\ap(M)=M$ and $C=C'M$, we may assume that $c\in C'$ and 
$\ap(c)\in CG^0=C'G^0$.   Then $\alpha(c)=c_1 g_1$, $c_1\in C'$, $g_1\in G^0$. Then
$\pi(\ap(c))=\pi(c_1)\pi(g_1)=\pi(g_1)\pi(c_1)$. Therefore, for all $n\in\N$,
$\pi(\ap(c^n))=\pi(\ap(c)^n)=\pi(c_1^n)\pi(g_1^n)=\pi(g_1^n)\pi(c_1^n)$. As $\pi(c_1^{-1}\ap(c))=\pi(g_1)$
generate a compact group in $\pi(U\ap(U)\cap CG^0)\cap \pi(G^0)\subset\pi(W\cap CG^0)\cap\pi(G^0)$. 
Therefore, $c_1^{-1}\ap(c)$ generate a compact group in $WM\cap G^0$. From the choice of $W$ as above, 
it follows that $c_1^{-1}\ap(c)\in M$.
As $M$ is contained in $C$ and $c_1\in C'\subset C$, $\ap(c)\in C$. \hfill{$\Box$}\eo

\bp\label{ct-td} Let $G$ be a locally compact group and let $\alpha\in
\Aut(G)$. Let $C_{G^0}(\ap)=\{x\in G\mid \ap^n(x)G^0\to G^0\mbox{ in } G/G^0\}$.
Let $M$ be an $\ap$-invariant compact subgroup of $G^0$ such that it is
normal in $G^0$ and $G^0/M$ is a Lie group.
Then the following hold:
\begin{enumerate}

\item[$(1)$] $M$ is normal in $C_{G^0}(\ap)$.
\item[$(2)$]  $C_{G^0}(\ap)=C_M(\ap)G^0$.
\item[$(3)$]  If $G^0$ is a Lie group, then $C_{G^0}(\ap)=C(\ap)G^0$.
\end{enumerate}
\ep

\bo Let $M$ be as in the hypothesis.
Let $V$ be a neighborhood of the identity $e$ in $G$.
Let $C\subset V$ be as in Lemma \ref{lc}.
Then $H=CG^0$ is an open subgroup in $G$, $C\cap G^0=M$, $H$ normalises $M$ and
$H/M=C/M\times G^0/M$. Let $x\in C_{G^0}(\ap)$.
There exists $N\in\N$ such that $\ap^n(x)\in H$ for all
$n\geq N$.  Let $y=\ap^N(x)$. Then
$y\in C_{G^0}(\ap)\cap H$ and $\ap^n(y)\in H$ for all $n\in \N$. In particular, $y=\ap^N(x)$ normalises $M$
and as $M$ is $\ap$-invariant, $x$ normalises $M$ and hence (1) holds.

For $y=\ap^N(x)$ as above, we have $y=cg$ for some $c\in C$ and $g\in G^0$ and
$\ap^n(y)=\ap^n(c)\ap^n(g)\in H$ for all $n\in\N$. Since $G^0$ is $\ap$-invariant, $\ap^n(c)\in H$.
Since $c\in C$, by Lemma \ref{lc}, $\ap^n(c)\in C$.
Moreover, as $\ap^n(y)G^0\to G^0$, we have that $\ap^n(c)G^0\to G^0$. As $H/M=C/M\times G^0/M$, we get that
$\ap^n(c)M\to M$ in $C/M$, i.e.\ $c\in C_M(\ap)\cap C$ and $y=cg\in (C_M(\ap)\cap C)G^0$. Now
$x=\alpha^{-N}(y)=\ap^{-N}(c)\ap^{-N}(g)\in C_M(\ap)G^0$, as both $C_M(\ap)$ and $G^0$ are $\ap$-invariant.
That is, $C_{G^0}(\ap)\subset C_M(\ap)G^0$. The converse is obvious as $M\subset G^0$. This proves (2).
Observe that  (2) $\Rightarrow$ (3) as one can take $M=\{e\}$ when $G^0$ is a Lie group.  \hfill{$\Box$}\eo

\bo $\!\!\!\!\!${\bf ~of Theorem \ref{clt-ctl}}
Let $G$, $\ap$ and $L$ be as in the hypothesis. We also know that
$C_{LG^0}(\ap)=C_{G^0}(\ap)L$ from the result on the totally
disconnected group $G/G^0$ (see Theorem 3.8 of \cite{BW} which is
valid for non-metrizable groups by \cite{J}). $C_L(\ap)\subset
C_{G^0}(\ap)L$ where $C_{G^0}(\ap)$ is an $\ap$-invariant subgroup
defined in Proposition \ref{ct-td}.

\medskip
\noindent{\bf Step 1} Now suppose $G^0$ is a Lie group. Then by
Proposition \ref{ct-td} (3), $C_{G^0}(\ap)=C(\ap)G^0$ and hence
$C_L(\ap)=C(\ap)(G^0\cap C_L(\ap))L$. Moreover, $G^0\cap C_L(\ap)=
G^0\cap C_{L\cap G^0}(\ap)=(C(\ap)\cap G^0)(L\cap G^0)$ (cf.\
\cite{HS}, Theorem 2.4). From the above equations we get that
$C_L(\ap)=C(\ap)L$. This shows that the assertion holds if $G^0$ is
a Lie group.

\medskip
\noindent{\bf Step 2} Since $C_L(\ap)\subset C_{G^0}(\ap)L$,
$C_L(\ap)=(C_L(\ap)\cap C_{G^0}(\ap ))L\subset (C_L(\ap)\cap
\ol{C_{G^0}(\ap )})L=(C_{L'}(\ap)\cap \ol{C_G^0(\ap)})L$, where
$L'=L\cap \ol{C_{G^0}(\ap)}$. Hence we may replace $L$ by $L'$ and
$G$ by $\ol{C_{G^0}(\ap)}$.

Suppose $G^0$ is $\ap$-Lie projective. Let ${\cal A} $ be the collection of all $\ap$-invariant 
compact normal subgroups $K$ of $G^0$ such that $G^0/K$ is a
Lie group.  Let $M \in {\cal A}$.  Then by Proposition \ref{ct-td} (1), $M$
is normal in $C_{G^0}(\ap)$ and hence in $G$. Applying Step 1 for
$G/M$, we get that $C_L(\ap)\subset C_{LM}(\ap)=C_{M}(\ap)L$.
As $G^0$ is $\ap$-Lie projective,
we have that $\cap {\cal A}=\{e\}$ and hence, since $\cal A$
is closed under finite intersection, we may assume that the family ${\cal A}$ is directed downwards. 
Now the assertion follows from Proposition 10 of \cite{J}.


\medskip
\noindent{\bf Step 3}  Suppose, the restriction of $\ap$
to the largest compact normal subgroup $K$ of $G^0$ is distal. Then by Lemma \ref{ind},
there exist $\ap$-invariant compact normal subgroups $K_i$ of $K$, such that
each $K/K_i$ is a Lie group and $\cap_i K_i=\{e\}$, (see also Theorem 2.6 and Corollary 2.7 in \cite{Jnew}).
As noted earlier in the proof of Proposition \ref{conn-ct}, $K_i$'s are normal in
$G^0$ and each $G^0/K_i$ is also a Lie group. Hence the assertion
holds from above. \hfill{$\Box$}\eo

\bc \label{cor2-4.2}
Let $G$ be a locally compact group and $\ap\in\Aut(G)$. Let $K$ be the largest compact normal subgroup of
$G^0$ and let $C$ be the largest closed (compact) $\ap$-invariant subgroup of $K$ such that the restriction of $\ap$ to
$C$ is ergodic. Then for any $\ap$-invariant subgroup $L$ of $G$, $C_L(\ap)=C_{L\cap C}(\ap)L$. In particular,
$C_L(\ap)=C(\ap)L$ if $L\cap C=\{e\}$.\ec

\bo As $K$ is characteristic in $G$, it is $\ap$-invariant and normal in $G$. Moreover, $(G/K)^0=G^0/K$ is a Lie group.
By Theorem 2.6 and Corollary 2.7 of \cite{Jnew}, $C$ is normal in
$K$ and it is also the smallest normal subgroup of $K$ such that the $\ap$-action on $K/C$ is distal. Moreover,
$K$ has a family ${\cal A}$
of compact normal $\ap$-invariant subgroups $M$ such that each $K/M$ is a Lie group and $\cap{\cal A}=C$; i.e.\ $K/C$ is
$\ap$-Lie projective. Since $G^0$ is connected and $K$ is a compact normal subgroup in it, we get that each $M$ in
${\cal A}$ is normal in $G^0$ and $G^0/M$ is a Lie group. As in the proof of Theorem \ref{clt-ctl}, we may replace $G$
by $\ol{C_{G^0}(\ap)}$ and $L$ by $L\cap\ol{C_{G^0}(\ap)}$ assume that every $M$ in ${\cal A}$, as well as $C$,
is normal in $G$ (cf.\ Proposition \ref{ct-td}) and $(G/C)^0=G^0/C$ is $\ap$-Lie projective. Let
$\bar\ap\in\Aut(G/C)$ be the automorphism corresponding to $\ap$. As $LC$ is a compact $\ap$-invariant subgroup,
using Theorem 4.2 for $G/C$,  we get that $C_{(LC/C)}(\bar\ap)=C(\bar\ap)(LC/C)$, and hence,
$C_{LC}(\ap)=C_C(\ap)L$. This implies that $C_L(\ap)=C_{L\cap C}(\ap)L$. In particular,
$C_L(\ap)=C(\ap)L$ if $L\cap C=\{e\}$. \hfill{$\Box$}\eo

\bo  $\!\!\!\!\!${\bf ~of Theorem \ref {ct}}: If the $\{\ap^n\}_{n\in\N}$-action on $G$ is distal, then $C(\ap ) =
\{ e \}$.  We now assume that $C(\ap ) =\{ e \}$.
Suppose $x\in G$ is such that $\ap^{n_k}(x)\to e$ for some unbounded
sequence $\{ n_k \}$. Let $H$ be the closure of $\{g\in G \mid \ap ^
{n_k}(g)\to e \}$. Then $H$ is a closed $\ap$-invariant subgroup and $x\in H$. Then
the restriction $\ap'$ of $\ap$ to $H$ is weakly contractive as it
has a dense subset $D$ such that $\ap ^{n_k}(g)\to e$ for all $g\in D$.
By Theorem 5 of \cite{J0}, there exists a compact
$\ap$-invariant subgroup $L$ of $H$ such that $H=C_L(\ap')$, i.e.
$H\subset C_L(\ap)$.  Since $C(\ap)$ is trivial, by Theorem 13.2 of \cite{Jnew},
the restriction of $\ap$ to any compact $\ap$-invariant subgroup is distal.  
By Theorem \ref{clt-ctl}, $C_L(\ap)=C(\ap)L=L$.  Therefore, 
$x\in L$ and $x=e$ as the restriction of $\ap$ to $L$ is distal. This shows that
the $\{ \ap^n \}_{n\in\N}$-action on $G$ is distal. \hfill{$\Box$}\eo

The following corollary follows easily from above and Theorem 1.1 of \cite{Sh}.

\bc \label{ct1} Let $G$ be a locally compact group and
$\ap\in\Aut(G)$. Then the closure of the $\ap$-orbit of $x$,
$\overline{\{\ap^n(x)\}}_{n\in\Z}$ is a minimal closed
$\ap$-invariant set for every $x\in G$ if and only if $C(\ap^{\pm
1})=\{e\}$. \ec

Note that for any locally compact group $G$ and the largest compact normal subgroup $K$ of $G^0$,
$G/K$ satisfies the hypothesis of Theorem \ref{clt-ctl} for any automorphism of it. Moreover, if $\ap$ is the
inner automorphism of $G$ by an element contained in an almost connected subgroup, then the
restriction of $\ap$ to $K$ as above is distal, and the decomposition theorem holds for such $\ap$. We now construct 
an example of a group $G$ with an automorphism $\ap$ which does not act distally on
the largest compact normal subgroup of $G^0$ such that $G$ is not Lie projective but $G^0$ is $\ap$-Lie projective, 
i.e.\ $G$ and $\ap$ satisfy the first condition in Theorem \ref{clt-ctl}. 

\begin{example}
Let $G=\Q_p\times (\Z\ltimes (N\times K))$, $p$ a fixed prime. Where $N$ is the 3-dimensional Heisenberg group of 
strictly upper triangular real matrices whose centre is the commutator subgroup
$N_1$, $K=(\T^2)^\Z$, where $\T^2$ is the 2-torus. Note that $N$ is a simply connected step-2 nilpotent group.
There is a canonical action of $SL_2(\Z)$ on $N$ and also on $\T^2$; the action on $N$ keeps elements of $N_1$ fixed.
The action of $\Z$ on $G^0=N\times K$ is defined as follows: For $z\in \Z$, $z(x, (t_n))z^{-1}=(\beta^z(x), (t_{n+z}))$ 
for all $x\in N$ and $(t_n)\in K$, where
$\beta=\begin{pmatrix} 2 & 3\\ 5 &7\end{pmatrix}\in SL(2,Z)$ and $z(t_n)z^{-1}=(t_{n+z})$ for every $(t_n)\in K$; it is the 
shift action of $\Z$ on $K$ which is ergodic. Let $\ap\in\Aut(G)$ be defined as follows:
$\ap(x)=px$ for all $x\in\Q_p$. For $z\in \Z$, $\ap(z)=z\gamma(z)$, where $\gamma$ identifies $\Z$ with a lattice of the 
centre $N_1$ of $N$, which is a discrete co-compact subgroup of $N_1$.
For $(x, (t_n))\in N\times K$, $\ap(x, (t_n))=(\beta(x), (\beta(t_n)))$. It is easy to see that $\ap$ is an automorphism of $G$. 
Note that $G$ is not Lie projective as $K$ does not have any nontrivial closed subgroup invariant under the action of 
$\Z$. But as $\ap$ keeps each $\T^2$ invariant, $K$ has $\ap$-invariant subgroups $K_i$ which are normal in $G^0$ 
such that each $G^0/K_i$ is a Lie group and $\cap_i K_i=\{e\}$, i.e.\ $G^0$ is $\ap$-Lie projective.
\end{example}

\end{section}

\begin{section}{Unimodularity}

In this section, we relate relative compactness of contraction
groups of inner automorphisms and groups whose closed subgroups are unimodular.

We first make the following observation that is an extension of
the well-known fact that any linear transformation of determinant
one preserves the Lebesgue measure.

\bp\label{pum} Any distal automorphism preserves the Haar measure.
In particular, any pointwise distal group is unimodular. \ep

\bo

Let $\ap$ be a distal automorphism of a locally compact group $G$ and let
$G^0$ be the connected component of $e$ in $G$.
Then by Theorem 3.3 of \cite{RS}, factor of $\ap$ on $G/G^0$ is also distal.
By Proposition 2.1 of \cite{JR}, $G/G^0$ has a $\ap$-invariant compact open subgroup
$K$ of  $G/G^0$ (see also \cite{BW} and \cite{J}).
Let $H$ be the closed subgroup of $G$ containing $G^0$ such that
$H/G^0=K$.  Then $\ap(H)=H$ and $H$ is an almost connected open subgroup.  Thus,
it is sufficient to prove that the Haar measure of $H$ is $\ap|_H$-invariant.
In other words, we may assume that $G$ is almost connected.  This implies
that there is a characteristic compact subgroup $M$ of $G$ such that
$G/M$ is a Lie group with finitely
many connected components. Therefore, $G^0M$ is an open $\ap$-invariant subgroup
of finite index in $G$.  Hence we may assume that $G=G^0M$. As  
$M$ is a compact normal $\ap$-invariant subgroup and $G^0M//M$ is a
connected Lie group, using Theorem 3.1 of \cite{RS}, we may assume that $G$ is a
connected Lie group.  Let $\du\ap$ be the differential of $\ap$ which
is a linear transformation on the Lie algebra of $G$.  Since $\ap$ is
distal on $G$, eigenvalues of $\du \ap$ have absolute
value one (see \cite{Ab1} and \cite{Ab2}).  This implies that
determinant of $\du \ap$ has absolute value one.  This proves that
$\ap$ preserves the Haar measure of $G$. \hfill{$\Box$}\eo

For a locally compact group $G$, we consider the following statements: 

\be

\item[$(i)$] Closed subgroups of $G$ are unimodular;

\item[$(ii)$] $C_L(\ap )$ is relatively compact for any inner automorphism
$\ap$ of $G$ and for any $\ap$-invariant compact subgroup $L$ of
$G$;

\item[$(iii)$] $C(\ap )$ is relatively compact for any inner
automorphism $\ap$ of $G$;

\item[$(iv)$] $G$ is pointwise distal.

\ee

We now obtain the following relations between $(i)$-$(iv)$.

\bt\label{csu} Let $G$ be a locally compact group.  Then $(i)$ and
$(ii)$ are equivalent and $(iv) \Rightarrow (ii) \Rightarrow (iii)$.

\be
\item[$(1)$]  Further, if the connected component in $G$ is a real Lie group, then
$(i)$--$(iii)$ are equivalent.

\item[$(2)$]  Furthermore, if there is a continuous injection $\vp \colon G\to GL (n, \mF
)$ where $\mF$ is a local field or if $G$ is a closed subgroup of an almost
connected group, then $(i)$--$(iv)$ are equivalent.
\ee \et

We first prove the following automorphism version of Theorem
\ref{csu}.

\bp\label{csua} Let $G$ be a locally compact group and let $\ap$ be
any automorphism of $G$. Then the following are equivalent:

\be

\item[$(i)$] The Haar measure of any $\ap$-invariant closed subgroup of
$G$ is $\ap$-invariant.

\item[$(ii)$] $C_L(\ap ^{\pm 1})$ is relatively compact for any
$\ap$-invariant compact subgroup $L$ of $G$.\ee \ep

\bo Suppose $C_L(\ap ^{\pm 1})$ is relatively compact for any
$\ap$-invariant compact subgroup $L$ of $G$.  To prove (i) it is
sufficient to prove that the Haar measure of $G$ is $\ap$-invariant.

\medskip
\noindent{\bf Case of real Lie group:}  We first assume that $G$ is
real lie group.  In this case $\ol {C(\ap ^{\pm 1} )}$ is
compact, hence they normalize each other. Moreover, 
they are connected subgroups of $G^0$. Let $K=\ol{C(\ap)C(\ap^{-1})}$. 
Then $K$ is a compact connected subgroup of $G^0$. 
Let $\cG$ be the Lie
algebra of $G^0$. Then $\du\ap$, the differential of $\ap$, defines a
linear transformation on $\cG$. Let $V\subset \cG$ be the Lie
algebra of $K$.  Then eigenvalues of the factor of $\du \ap$ on $\cG
/V$ have absolute value one and hence its determinant has absolute
value one. Since $V$ is a Lie algebra of a compact connected Lie
group, the determinant of $\du \ap$ restricted to $V$ also has
absolute value one. Thus, the absolute value of determinant of $\du
\ap$ on $\cG$ is one.  This proves that the Haar measure on $G$ is
$\ap$-invariant.

\medskip
\noindent {\bf General case:}  Let $G^0$ be the connected component
of the identity $e$ in $G$. Then $G^0$ is $\ap$-invariant and $G/G^0$ is totally
disconnected.  We denote the factor of $\ap$ on $G/G^0$ by $\ba$
which is an automorphism defined by $\ba (gG^0) = \ap (g) G^0$ for
all $g\in G$. Let $\pi \colon G \to G/G^0$ be the canonical
projection.  Then $\pi (\ap (g) ) = \ba (\pi (g))$ for all $g\in G$.
Let $s_{G/G^0} \colon \Aut (G/G^0) \to \N$ be the scale function
defined as in \cite{BW}. Proposition \ref{ct-td} (2) together with the hypothesis
implies that $C(\ba^{\pm 1})$ is relatively compact in $G/G^0$, hence by
Proposition 3.24 of \cite{BW} we get that $s_{G/G^0}(\ba ^{\pm 1}) =1$.  This
implies that there is a $\ap$-invariant open subgroup $H$ of $G$
such that $G^0\subset H$ and $H/G^0$ is compact. Therefore, $H$ is
almost connected.  By \cite{MZ}, $H$ has a maximal compact normal
subgroup $K$ such that $H/K$ is a real Lie group. Since $K$ is
maximal, it is characteristic.  In particular, $K$ is
$\ap$-invariant.  Let $\nu$ be the factor automorphism of $\ap$ on
$H/K$. Then $C_K(\ap )/K = C(\nu )$ and hence $C(\nu )$ as well as
$C(\nu ^{-1})$ is relatively compact.  Since $H/K$ is a real Lie
group, we get from the real Lie group case that the Haar measure of
$H/K$ is $\nu$-invariant. As $K$ is compact, the Haar measure of
$H$ is $\ap$-invariant. Since $H$ is an open subgroup of $G$, the
Haar measure of $G$ is $\ap$-invariant.

Conversely, assume that the Haar measure of any $\ap$-invariant
closed subgroup of $G$ is $\ap$-invariant.  Let $L$ be any
$\ap$-invariant compact subgroup of $G$.  Then by Theorem 5 of
\cite{J0}, $\ol {C_L(\ap )} = C_{L_1}(\ap )$ for some
$\ap$-invariant compact subgroup $L_1$ of $G$. By assumption, Haar
measure of $C_{L_1}(\ap )$ is $\ap$-invariant, hence by Proposition
3 of \cite{J0}, we get that $C_{L_1}(\ap )= L_1$.  This proves that
$C_L(\ap )$ is relatively compact. \hfill{$\Box$}\eo

In order to complete the proof of Theorem \ref{csu}, we need the
following lemma.

\bl\label{inr} Assume that closed subgroups of $G$ are unimodular.
Then for any inner automorphism $\ap$ and any $\ap$-invariant closed
subgroup $H$ of $G$, the Haar measure on $H$ is $\ap$-invariant. \el

\bo Let $\ap$ be an inner automorphism of $G$ defined by $x\in G$.
Suppose $H$ is any $\ap$-invariant closed subgroup of $G$.  Then let
$H_1$ be the closed subgroup generated by $x$ and $H$.  Since $H$ is
$\ap$-invariant, $H$ is normal in $H_1$. Since $H_1/H$ is the closed
subgroup generated by a single element, $H_1/H$ is either discrete
or compact.

If $H_1/H$ is discrete, then $H$ is open in $H_1$.  By assumption
$H_1$ is unimodular and since $x\in H_1$, the Haar measure on $H_1$ is
$\ap$-invariant.  Since $H$ is open in $H_1$, the Haar measure on $H$ is
also $\ap$-invariant.

If $H_1/H$ is compact, let $\lam _H$ denote the Haar measure on $H$
and $\delta : H _1\to (0, \infty )$ be such that $\lam _H(gEg^{-1})=
\delta (g)\lam _H (E)$ for any Borel subset $E$ of $H$.  Then
$\delta $ is a continuous homomorphism.  By assumption $H$ is
unimodular, hence $\delta (H)= \{ 1 \}$.  Since $H_1/H$ is compact,
$\delta (H_1)$ is a compact subgroup of $(0, \infty )$, hence
$\delta (H _1)= \{ 1 \}$. \hfill{$\Box$}\eo

\bo $\!\!\!\!\!${\bf ~of Theorem \ref{csu}~} Applying Proposition
\ref{csua} to inner automorphisms, we get that $(ii)$ implies $(i)$,
and $(i)$ implies $(ii)$ follows from Lemma \ref{inr} and by
applying Proposition \ref{csua} to inner automorphisms.

By Corollary 3.2 of \cite{RS}, we get that $(iv)$ implies $(ii)$ and
by taking $L= \{ e\}$ we get that $(ii)$ implies $(iii)$.

If the connected component in $G$ is a real Lie group.  Then $(iii)
\Ra (ii)$ follows from Theorem \ref{clt-ctl}.

Suppose there is a continuous injection $\vp \colon G\to GL (n, \mF
)$ where $\mF$ is a local field. In order to prove $(i)$--$(iv)$ are
equivalent it is sufficient to prove that $(iii)$ implies $(iv)$.
Let $x\in G$ and define $\ap \colon GL (n, \mF )\to GL (n, \mF )$ by
$\ap (g) = \vp (x) g \vp (x^{-1})$ for all $g\in GL (n, \mF )$. Then
$\vp (xgx^{-1}) = \vp (x) \vp (g) \vp (x^{-1})$ for all $g\in G$.
Let $C(\ap , GL (n, \mF )) = \{ g \in GL (n, \mF ) \mid \ap ^n (g)
\to e ~~{\rm as}~~ n \to \infty \}$ and $C(\ap , G) = \{ g \in G
\mid x ^n g x^{-n}\to e ~~{\rm as}~~ n \to \infty \}$.  Then $C(\ap
, GL (n, \mF ) )$ is a closed $\ap$-invariant subgroup of $GL (n,
\mF )$ and $\vp (C(\ap , G)) \subset C(\ap, GL (n, \mF ))$.  Since
$\vp$ is continuous, $\vp (\ol {C(\ap , G)})\subset \ol {\phi (C(\ap
, G))} \subset C(\ap , GL (n, \mF ))$.  Now, $(iii)$ implies that
$\ol {C(\ap , G)}$ is a compact subgroup of $G$, hence $\vp (\ol
{C(\ap , G)})$ is a $\ap$-invariant compact subgroup.  Since $C(\ap,
GL(n, \mF))$ has no $\ap$-invariant nontrivial compact subgroup,
$\vp (\ol{C( \ap, G)})$ is trivial, hence $C(\ap, G)$ is trivial as
$\vp$ is an injection.  Now $(iv)$ follows from Theorem \ref{ct}.

Suppose $G$ is a closed subgroup of an almost connected group $H$.
Then $H$, and hence, $G$ is Lie projective, i.e.\ $H$ (resp.\ $G$) has arbitrarily
small compact normal subgroups $K'_d$ (resp.\ $K_d$) such that $H/K'_d$ (resp.\
$G/K_d$)  is a real Lie group and $H/K'_d$ has finitely many connected components.
By Theorem \ref{clt-ctl}, each $G/K_d$ satisfies $(iii)$ if $G$ satisfies $(iii)$.
Moreover, if each $G/K_d$ is pointwise distal, then so is $G$.
Thus, it is sufficient to prove $(iii)$ implies $(iv)$ when $G$ is contained in a
real Lie group $H$ such that $H$ has finitely many connected components.
Now as $G\cap H^0$ is a closed normal subgroup of finite index in $G$, to
prove $(iii)$ implies $(iv)$, we may further assume that $G$ is contained in
a connected Lie group $H$. Suppose $G$ satisfies $(iii)$. Let $x\in G$. Define $\ap \colon H \to H$ by
$\ap (g)= xgx^{-1}$ for all $g \in H$. As $G$ satisfies $(iii)$, $C(\ap )\cap G$ and $C(\ap
^{-1})\cap G$ are relatively compact. We will now claim that $\ap|_G$ is
distal.  Let $K$ be the maximal compact normal subgroup of $H$.
As $H$ is connected, the action of $H$, and hence, the action of $G$ on $K$ by
conjugation is by inner automorphism of $K$ (cf.\ \cite{I}, Theorem $1'$)
and hence, in particular, the
action of the inner automorphism $\ap$ of $H$ on $K$ is distal.
Then by Proposition \ref{conn-ct}, $C(\ap^{\pm 1})$ is closed, and hence,
a simply connected nilpotent Lie subgroup of $H$, in particular, they do
not have any non-trivial compact subgroups. Moreover, $C(\ap)\cap G$ and
$C(\ap^{-1})\cap G$ are closed, and hence, compact subgroups of $C(\ap)$ and
$C(\ap^{-1})$ respectively, and hence they are trivial.
Now Theorem \ref{ct} implies that $\ap|_G$ is distal (see also \cite{Ab1}, \cite{Ab2}).
As this is true for any $x\in G$, we have that $G$ is pointwise distal.
Thus, $(iii)$ implies $(iv)$.
\hfill{$\Box$}\eo

\end{section}

\begin{acknowledgement} The second author would like to thank the Indian Statistical Institute (ISI),
Bangalore, India as some part of the work was done during a couple of short visits to the institute. Moreover,
she would also like to thank the National Board for Higher Mathematics (NBHM), DAE, Govt.\ of India for a project
grant during the period when the manuscript was revised and improved.
\end{acknowledgement}

\bigskip\medskip
\advance\baselineskip by 2pt
\begin{tabular}{ll}
C.\ R.\ E.\ Raja & \hspace*{1cm}Riddhi Shah \\
Stat-Math Unit & \hspace*{1cm}School of Physical Sciences(SPS)\\
Indian Statistical Institute (ISI) & \hspace*{1cm}Jawaharlal Nehru University(JNU)\\
8th Mile Mysore Road & \hspace*{1cm}New Delhi 110 067, India\\
Bangalore 560 059, India & \hspace*{1cm}rshah@jnu.ac.in\\
creraja@isibang.ac.in & \hspace*{1cm}riddhi.kausti@gmail.com
\end{tabular}

\end{document}